\documentclass{itor}

\usepackage[authoryear]{natbib}%
\usepackage[figuresright]{rotating}

\usepackage{amsmath,amsthm}
\newtheorem{theorem}{Theorem}

\theoremstyle{definition}

\theoremstyle{remark}
\newtheorem{definition}{Definition}

\usepackage{tikz}
\usetikzlibrary{automata}
\usepackage[short,nocomma]{optidef}
\usepackage[capitalise]{cleveref}
\usepackage[linesnumbered]{algorithm2e}
\usepackage{subcaption}
\usepackage{amsfonts}
\usepackage{float}
\usepackage{placeins}
\usepackage{multirow}
\usepackage{multicol}
\usepackage{adjustbox}
\begin{document}

% Title
\title{The Agricultural Spraying Vehicle Routing Problem With Splittable Edge Demands}

%Authors, affiliations address.
\author[Running Author]{Qian Wan\affmark{a,b$\ast$}, Rodolfo Garc\'ia-Flores \affmark{b}, Simon A. Bowly \affmark{c}, Philip Kilby \affmark{d} and Andreas T. Ernst \affmark{a}}

\affil{\affmark{a}School of Mathematics/Monash University/Wellington Rd, Clayton, 3800, VIC, Australia.}
\affil{\affmark{b}Data Driven Enterprises/CSIRO Data-61/Gate 5 Normanby Rd, Clayton, 3168, VIC, Australia.}
\affil{\affmark{c}Gurobi Optimization/9450 SW Gemini Dr. \#90729, Beaverton, 97008-7105, Oregon, USA.}
\affil{\affmark{d}Data61/CSIRO/GPO Box 1700, Canberra, 2601, ACT, Australia.}
\email{Qian.Wan@momash.edu;\\ Rodolfo.Garcia-Flores@data61.csiro.au;\\ bowly@gurobi.com; \\Philip.Kilby@data61.csiro.au;\\Andreas.Ernst@monash.edu. }

%Author's notes/correspondence etc.
\thanks{\affmark{$\ast$}Author to whom all correspondence should be addressed (e-mail: Qian.Wan@momash.edu).}

%%Date
\historydate{Received DD MMMM YYYY; received in revised form DD MMMM YYYY; accepted DD MMMM YYYY}

%Abstract
\begin{abstract}
In horticulture, spraying applications occur multiple times throughout any crop year. This paper presents a splittable agricultural chemical sprayed vehicle routing problem and formulates it as a mixed integer linear program. The main difference from the classical capacitated arc routing problem (CARP) is that our problem allows us to split the demand on a single demand edge amongst robotics sprayers. We are using theoretical insights about the optimal solution structure to improve the formulation and provide two different formulations of the splittable capacitated arc routing problem (SCARP), a basic spray formulation and a large edge demands formulation for large edge demands problems. This study presents solution methods consisting of lazy constraints, symmetry elimination constraints, and a heuristic repair method. Computational experiments on a set of valuable data based on the properties of real-world agricultural orchard fields reveal that the proposed methods can solve the SCARP with different properties. We also report computational results on classical benchmark sets from previous CARP literature. The tested results indicated that the SCARP model can provide cheaper solutions in some instances when compared with the classical CARP literature. Besides, the heuristic repair method significantly improves the quality of the solution by decreasing the upper bound when solving large-scale problems.
\end{abstract}

%Keywords, etc.

\maketitle

% Heading 1
\section{Introduction}\label{sec1}

% Heading 2
\subsection{Motivation}
    \noindent Precision agriculture has been implemented to provide low-input, high-efficiency, and sustainable agricultural production. In precision agriculture, automation and robotics have become common practice in planting \citep{cay2018development,khazimov2018development,shi2019numerical}, inspection \citep{mahajan2015image,singh2017detection,ozguven2019automatic}, spraying \citep{oberti2016selective,vakilian2017farmer,gonzalez2016autonomous}, and harvesting \citep{xiong2019development,williams2019robotic,barth2016design} to overcome the labour shortcomings and improve efficiency \citep{mahmud2020robotics}.
        
    Spraying applications are one of the most critical sectors of agricultural areas. Recently, autonomous agricultural robotic sprayers have been employed in agricultural lands to perform spraying applications such as spraying pesticides, herbicides, and fungicides. Robotic sprayers are a sustainable and human-friendly way to perform spraying operations in orchards. Moreover, as the scale of the orchard increases, the required time frame for spraying operations remains the same. As a result, it is unimaginable to perform the spaying applications only by human workers in a limited time for a large size orchard. Therefore, robotic sprayers have become a more and more popular way to achieve spraying operations in orchard areas.
        
    Finding an optimal routing plan for robotic sprayers is one of the main work in spraying operations. Optimizing the route plan of a robotic sprayer has many benefits, such as reducing time and overall cost, requiring fewer workers, and thus increasing profits \citep{xu2022efficiency}.
        
    This paper presents a mixed integer programming-based optimization model for an orchard robotic sprayer routing problem. We give two formulations for the optimization model: the first formulation is a basic spray formulation, and the second is a large edge demands formulation. Our problem allows for multiple robots to satisfy the demand of a row by splitting the spray load. This leads to more efficient utilisation of the tank capacity of the robotic sprayer. We will show that this can provide a cheaper solution than the case where sprayers are required to always process complete rows. We aimed to optimize the path taken by the automatic sprayer to reduce the number of chemicals sprayed while maintaining their effectiveness.
    
% Heading 3 
\subsection{Problem definition}          
    \noindent The classical capacitated arc routing problems (CARP) aim to find minimum-cost paths for a fleet of identical vehicles of capacity P that must service the demand of a subset of edges in a graph. The CARP consists of finding a set of vehicle trips that must satisfy: 1) each vehicle departs from the depot, services a subset of the required edges, and returns to the depot; 2) the total demand serviced by a vehicle does not exceed P; 3) each required edge is serviced by exactly one vehicle. The CARP model is a good starting point as an abstraction of our problem, where vehicles are spraying robots and the arcs to be serviced are rows of an orchard, with demand of an arc being the amount of liquid that has to be sprayed in that row. However, the requirement that only a single vehicle must serve a whole arc seems unnecessarily restrictive. This paper solves an extension of CARP, which allows splitting the demand on a single demand edge amongst robots.
        
\begin{figure}[htb!]
            \centering
            \includegraphics[scale=0.5]{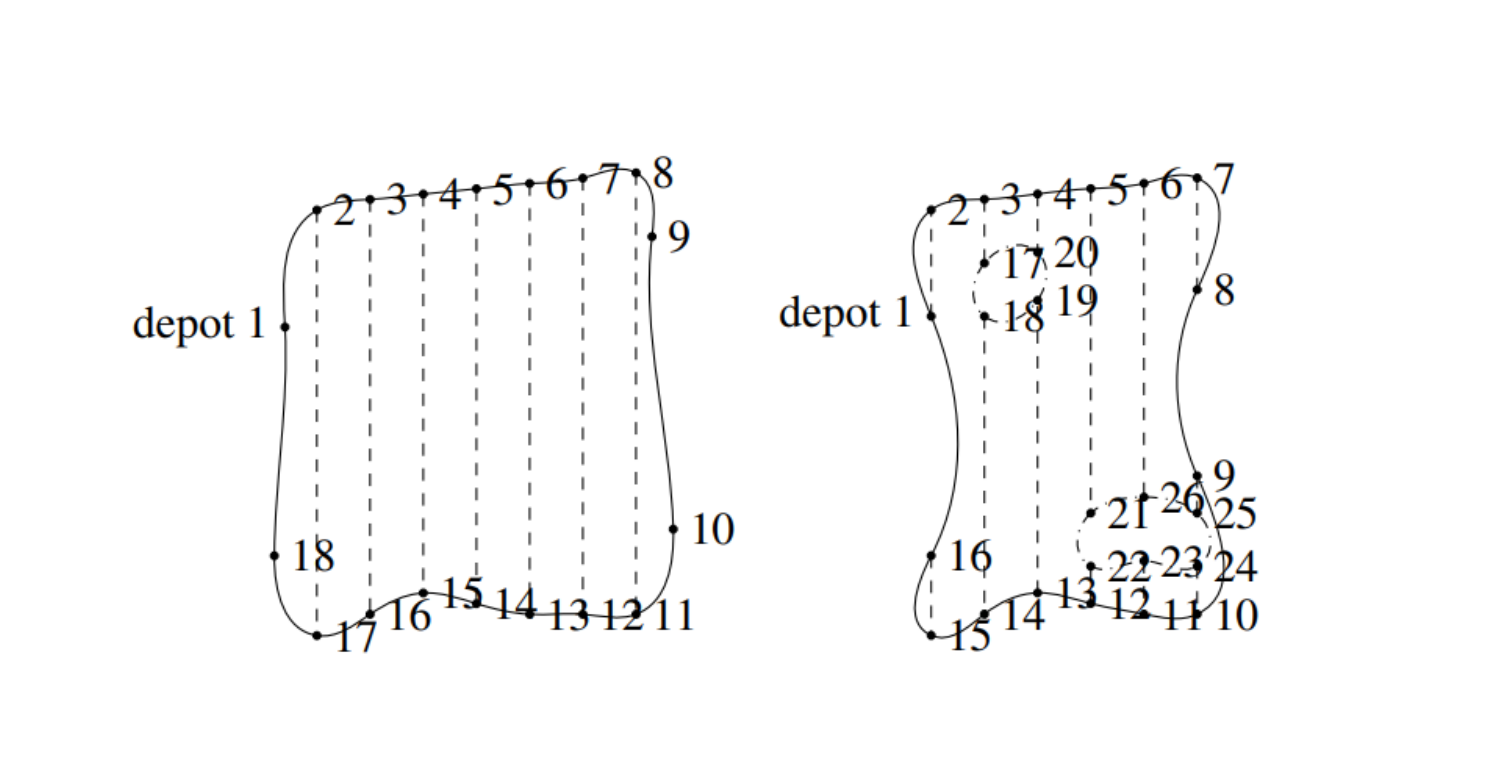}
        \caption{ An example of the small orchards based on the attributes of the real-world examples \citep{plessen2019optimal}: solid, dashed, and dash-dotted lines denote the headland, interior, and island edges. Dots indicate the vertices. (Left) Uninterrupted edges when aligned in a rotated coordinated frame. (Right) Interrupted edges due to orchard field indentations and obstacle areas prohibit the passage of any vehicle operating in the orchard. We set vertex 1 as the \emph{depot}, and robotic sprayers depart from the depot and spray the interior edges (dashed line). If the robotic sprayers run out of chemicals, it goes back to the depot to refill.}
        \label{Fig.1}
    \end{figure}

    This paper presents a splittable capacitated arc routing problem (SCARP) for orchards. The real-world orchard's shape, tractor tracks, barrier blocks, and tracks connecting the barriers form the graph of the orchard. An orchard is represented as an undirected graph $\mathcal{G}$. The tractor tracks are the lines on the graph $\mathcal{G}$. In graph theory, we call these lines as \emph{edges}. In practice, some \emph{edges} that need to be sprayed are \emph{demand edges} or \emph{required edges}, with each required edge having a demand, and those without demanding are \emph{non-demand edges}. Edges in the orchard were classified as headland edges, interior edges, and island edges. In practice, we only need to spray interior edges (for example, dashed lines in \cref{Fig.1}). Any intersection between two edges on a graph is \emph{vertex}. In practice, huge barrier blocks are formed by planted trees and islands that prohibit vehicles from passing. We call these barriers as \emph{obstacle areas}. \cref{Fig.1} is an example of the small orchards based on the attributes of the real-world examples based on the graphs from ~\cite{plessen2019optimal}. \cref{Fig.1} provides two main types of orchard graphs: uninterrupted (no obstacles) and interrupted graphs. The dots are vertices, whereas the solid, dashed, and dash-dotted lines represent the headland, interior, and island edges. In this problem, all robots start from a fixed depot (for example, vertex $1$ is a \emph{depot} in \cref{Fig.1}). The robot does spraying work subject to the robot's capacity and edge demands of this orchard. Every required edge must be sprayed by at least one robotic sprayer. The robotic sprayers are refilled at the depot when they run out of chemicals to spray. It is worth noting that this problem could be seen as both a multi-robotic and single-robotic spraying routing problem.
        
    In summary, the properties of the graph and robot sprayer of SCARP are:
    \begin{enumerate}
        \item Arbitrarily shaped orchards (convex or non-convex)
        \item Availability of a transition graph with edges connecting vertices
        \item A corresponding  non-negative cost array with edge costs equal to the path length of edges
        \item Corresponding non-negative edge demands (Edge demand should be non-zero for those edges that must be sprayed)
        \item A fixed robotics capacity. The total demand serviced by a robot does not exceed the capacity
        \item A designated start vertex (depot), same start and end vertex
        \item Fix the depot (refilled position) at the start vertex, and if the robot runs out of liquid to spray, it goes back to refill
        \item The demand of one edge may be satisfied by more than one robot (split demand)
        \item There is a single robotic sprayer by default
        \end{enumerate}
       
    The SCARP is an \emph{NP-hard} problem, since it can be reduced to the \emph{NP-hard} Rural Postman Problem (RPP) \citep{lenstra1976general} by giving the robot an unlimited capacity. This study presents solution methods consisting of lazy constraints, symmetry elimination constraints, and a heuristic repair method. We tested our algorithm on valuable data sets based on a real-world problem and a set of classical instances from CARP literature. This research proposes several contributions to developing an agricultural robot performing spraying tasks autonomously and gives some inspiration for mobile routing problems.

    The remaining paper is organized as follows: A brief review of the relevant literature is given in the~\cref{sec2}. Next, we provide a mathematical formulation of the SCARP in~\cref{sec3}. After that, we analyze the properties of the SCARP in~\cref{sec4}. The solution strategy is given in~\cref{sec5}. In addition, we reveal the success of the algorithms in~\cref{sec6}. Then, \cref{sec7} shows the parameter analysis of the heuristic repair method. Finally, we conclude the paper and point out the future research directions in~\cref{sec8}.
    
% Heading 1        
\section{Literature review}\label{sec2}
    Optimizing logistics and route planning plays a vital role in the agricultural area. The minimum navigation cost by finding the shortest path in mobile robot navigation could provide a convenient service with the lowest cost to the farmers. 
            
    A capacitated variation of Vehicle Routing Problem (VRP) is known as the Capacitated Vehicle Routing Problem (CVRP) \citep{braekers2016vehicle}. The CVRP seeks total cost-minimizing routes for multiple identical vehicles with limited capacity and where multiple vertices (customers) subject to various demands must be serviced exactly once by exactly one vehicle \citep{laporte1992vehicle}. \cite{archetti2012vehicle} propose the Split Delivery Vehicle Routing Problem (SDVRP), where each customer may be served by multiple vehicles. 
    The VRPs are computationally expensive and categorized as NP-hard \citep{lenstra1981complexity}. There are many variations \citep{lysgaard2020matheuristic,chiang2014knowledge,toth2014vehicle,laporte1992vehicle}. The well-known Traveling Salesman Problem (TSP) \citep{grotschel1985polyhedral} is a particular case of the VRP when only one vehicle is available at the depot, and no additional operational constraints are imposed. Numerous applications exist that involve integrating VRP with other problems, such as combining it with the depot location problem to construct a Location Routing Problem (LRP) \citep{quintero2017biased}. In contrast to the robotics spraying agricultural routing problem focus on \emph{edge-coverage}, the VRPs typically \emph{vertex-coverage} is of primary importance. Therefore, instead of VRPs, arc routing problems (ARPs) are here of more interest. 
            
    In Arc Routing Problem (ARP), the customers are located not at the vertices but along the arcs of the road network, are described in the book edited by \cite{dror2012arc}. There are many real-life applications of ARP such as garbage collection \citep{amponsah2004investigation}, postal delivery \citep{eiselt1995arc_a, eiselt1995arc_b}, snow removal \citep{campbell2000roadway} and others. \cite{monroy2017rescheduling} introduce a dynamic ARP that addresses disruptions caused by vehicle failure in the incapacitated case. For ARPs, problems can be distinguished between the classes of the Chinese Postman Problem (CPP) and the Rural Postman Problem (RPP), where all and only a subset of all arcs of the graph need to be traversed, respectively \citep{eiselt1995arc_a,eiselt1995arc_b}. However, neither the CPP nor the RPP considers the vehicle's capacity or the edge demands. 
            
    The Capacitated Arc Routing Problem (CARP) was formally introduced by \cite{golden1981capacitated}. This important problem plays the same role in arc routing as the CVRP in node routing. Roughly speaking, the CARP consists of finding a set of minimum-cost routes for vehicles with limited capacity that service some street segments with known demand, represented by a subset of edges of a graph. The CARP is an NP-hard problem. \cite{golden1981capacitated} showed that even the 0.5 approximate CARP is NP-hard. Besides, finding a feasible solution for the CARP that uses a given number of vehicles is also NP-hard. Thus, the CARP is considered to be very difficult to solve exactly. The SCARP can be seen as a relaxation problem of the classical CARP. The difference between the CARP and the SCARP in this paper is that, in SCARP, the demand for a single edge may be satisfied by more than one vehicle. 
            
   \cite{benavent1992capacitated} developed a dynamic programming-based technique to compute lower bounds to the CARP with a fixed number of vehicles and generated a set of network instances for the CARP. Our research also works on a fixed precomputed number of vehicles, and we provide a comparison result with these Benavent's instances set. Some authors have solved the CARP by transforming it into a capacitated node routing problem \citep{baldacci2006exact}, which applies an exact algorithm to compute the lower bounds of the CARP. \cite{belenguer2015chapter} classified exact methods into three groups: branch-and-bound based on combinatorial lower bounds, cutting-plane and branch-and-cut, column-generation and branch-and-price methods. \cite{belenguer2003cutting} introduced some new valid inequalities for the CARP and implemented a cutting plane algorithm for this problem. However, increasing the lower bounds significantly for large-scale problems by exact methods is difficult. Moreover, exact methods are limited to moderated instances \citep{chen2016hybrid}. Several heuristic algorithms have been developed for the CARP, such as simulated annealing algorithm \citep{eglese1994routeing}, tabu search heuristic \citep{hertz2000tabu}, improved ant colony optimization based algorithm \citep{santos2010improved}, etc. \cite{lacomme2004competitive} proposed memetic algorithms (MAs), which provided a better solution in some classical instance sets for solving an extended CARP. This paper provides a heuristic repair method to decrease the upper bounds for SCARP, which could provide cheaper solutions in some large-scale CARP instances.
            
    In this paper, we present a SCARP problem and formulate it as a mixed integer linear program. We are using theoretical insights about the optimal solution structure to improve the solutions. We implement lazy constraints, symmetry elimination constraints, and heuristic repair method to solve the SCARP. Computational experiments show that the SCARP model can provide cheaper solutions in some instances when compared with the classical CARP literature. Besides, the heuristic repair method significantly improves the quality of the solution by decreasing the upper bound when solving large-scale problems. 

\section{Mathematical formulations}\label{sec3}

\subsection{Nomenclature}
    \begin{itemize}
        \item Sets
            \begin{itemize}
                \item[] $\mathcal{V}=\{ 1,2,...,n\}$ set of all vertices;
                \item[] $\mathcal{E}=\{(i,j):i,j\in \mathcal{V},\ i\neq j\}$ set of all existing directed arcs of this graph. When the distance from $i$ to $j$ is not zero, it means this arc exists;
                \item[] $\mathcal{U}=\{(i,j)\in \mathcal{E} :i < j\}$  set of undirected arcs (set of arcs, when not need to consider the direction). 
                \item[] $\mathcal{G=(V,U)}$ an undirected connected graph;
                \item[] $\mathcal{I}_k=\{i\in \mathcal{E} :(i,k)\in \mathcal{E}\}$ set of incoming arcs of vertex $k$;
                \item[] $\mathcal{Q}_k=\{j\in \mathcal{E} : (k,j)\in \mathcal{E}\}$ set of outgoing arcs of vertex $k$;
                \item[] $\mathcal{R}=\{ 1,2,...,r\}$  set of the index of the \emph{robotic sprayer's tour} for single robotic sprayer problem or \emph{robotic sprayer} for multiple robotic sprayers problem. $r$ equals the smallest integer greater than the total demand divide the capacity of the robot sprayer $r=\frac{\sum_{(i,j)\in \mathcal{V}}D_{i j}}{P}$. We can refer to it as \emph{robot} to simplify.
                 % \item $\mathcal{A}=\{(i,j):i,j\in N,\ i\neq j\}$ Set of directed arcs connecting vertex $i$ to vertex $j$,  $(i,j) \neq (j, i)$;
            \end{itemize}
        \item Parameters
            \begin{itemize}
                \item[] $C_{i j}$ non-negative cost (km) of arc $(i,j)\in \mathcal{U}$;
                \item[] $D_{i j}$ non-negative demand of spray (L) of arc $(i,j)\in \mathcal{U}$;
                \item[] $P$ robot sprayer capacity (L);
                \item[] fixed depot is \emph{vertex 1} by default.
            \end{itemize}
        \item Decision variables
            \begin{itemize}
                \item[] $x_{i j}^r$ a binary variable, which is 1 if the arc $(i,j)\in \mathcal{E}$ is traversed by $r \in \mathcal{R}$, 0 otherwise;
                \item[] $y_{i j}^r$ the amount of chemicals sprayed of edge $(i,j)\in \mathcal{U}$ with robot  $r\in \mathcal{R}$. Fix $y_{i j}^r = 0$ if $D_{i j} = 0$.
            \end{itemize}
        \end{itemize}
        
\subsection{SCARP formulation 1: Basic robotic spraying formulation}
    \begin{mini!}|s|[1]
        {}{\sum_{r\in\mathcal{R}}\sum_{(i,j)\in\mathcal{U}}C_{i j} (x_{i j}^r+x_{j i}^r) \quad  (j,i)\in\mathcal{E} }{}{\label{Const0}}
        \addConstraint{\sum_{(i,j)\in\mathcal{U}} y_{i j}^r}{\leq P}{\quad \forall\ r\in \mathcal{R}\label{Const1}}
        \addConstraint{y_{i j}^r}{\geq 0}{\quad \forall\ (i,j)\in\mathcal{U},\ r\in\mathcal{R}\label{Const2}}
        \addConstraint{\sum_{r \in \mathcal{R}}y_{i j}^r}{=D_{i j}}{\quad\forall\ (i,j)\in \mathcal{U}\label{Const3}}
        \addConstraint{y_{i j}^r}{\leq P(x_{i j}^r+x_{j i}^r)}{\quad\forall\ (i,j)\in \mathcal{U},\ r \in \mathcal{R}\label{Const4}}
        \addConstraint{\sum_{i\in\mathcal{I}_k}x_{i k}^r}{=\sum_{j\in \mathcal{Q}_k}x_{k j}^r}{\quad\forall\ k\in \mathcal{V},\ r\in \mathcal{R}\label{Const5}}
        \addConstraint{\sum_{i\in \mathcal{S},j\in\mathcal{T}}Px_{i j}^r}{\geq \sum_{i,j\in \mathcal{T}}y_{i j}^r}{\ \ \begin{array}{l@{}l}\forall&\ \mathcal{S}\subset \mathcal{V},\ 1\in\mathcal{S},
            %\mathcal S\cup\mathcal{T}=\mathcal{N},%AE: this is implied by T=N\S
            \\&\mathcal{T}=\mathcal{V}\setminus\mathcal{S},\ r \in\mathcal{R}\end{array}\label{Const6}}
        \addConstraint{ x_{i j}^r}{\in\{0,1\}}{\quad\forall\ (i,j)\in \mathcal{E},\ r \in \mathcal{R}\label{Const7}}
    \end{mini!}

    The objective function (\ref{Const0}) seeks to minimize the total distance travelled. Constraints (\ref{Const1}) state the total sprays of a robot do not exceed the robot's capacity. Constraints (\ref{Const2}) confirm the number of sprays on each edge are greater than or equal to 0. Constraints (\ref{Const3}) ensure the total amount of sprays on each edge equal the corresponding edge demands. Constraints (\ref{Const4}) guarantee that the edge $(i, j)$ can be sprayed by a robot only if the robot covers edge $(i, j)$ with robot $r$ (robot $r$ can serve edge  $(i, j)$ in both directions). Constraints (\ref{Const5}) are flow conservation constraints for each vertex, which assure route continuity. Constraints (\ref{Const6}) enforce that any set of vertices $\mathcal T$ that contains arcs are sprayed by a robot $r$ ($y_{ij}^r>0$) must be reachable by the robot from the other side of the cut that contains the depot (set $\mathcal S$). Integrality restrictions are given in Constraints (\ref{Const7}).

\section{Properties of the SCARP}\label{sec4}

\subsection{Properties}
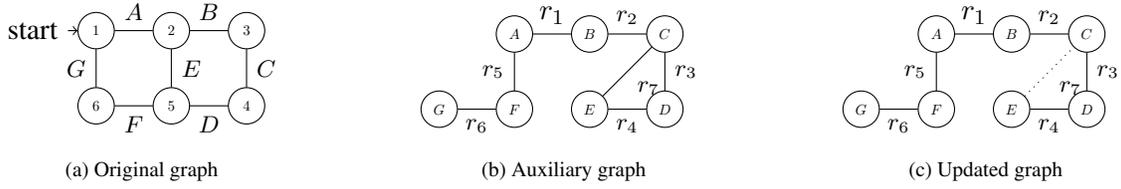
\begin{figure}
        \centering
        \begin{subfigure}{0.3\textwidth}
            \centering
            \begin{tikzpicture}[-,node distance=2cm,auto,scale = 0.25]        
            \node[initial,state, scale=0.5](A){$1$};
            \node[state,scale=0.5](B)[right of=A]{$2$};
            \node[state,scale=0.5](C)[right of=B]{$3$};
            \node[state,scale=0.5](D)[below of=C]{$4$};
            \node[state,scale=0.5](E)[below of=B]{$5$};
            \node[state,scale=0.5](F)[below of=A]{$6$};
            \path (A)edge node{\footnotesize$A$}(B)
            (B)edge node{\footnotesize$B$}(C)
            (C)edge node{\footnotesize$C$}(D)
            (D)edge node{\footnotesize$D$}(E)
            (B)edge node{\footnotesize$E$}(E)
            (E)edge node{\footnotesize$F$}(F)
            (F)edge node{\footnotesize$G$}(A);
            \end{tikzpicture}
            \caption{Original graph}
         \label{Original graph}
        \end{subfigure}
     \hfill
     \begin{subfigure}{0.3\textwidth}
         \centering
         \begin{tikzpicture}[-,node distance=2cm,auto,scale = 0.25]
        \node[state,scale=0.5](A){$A$};
        \node[state,scale=0.5](B)[right of=A]{$B$};
        \node[state,scale=0.5](C)[right of=B]{$C$};
        \node[state,scale=0.5](D)[below of=C]{$D$};
        \node[state,scale=0.5](E)[below of=B]{$E$};
        \node[state,scale=0.5](F)[below of=A]{$F$};
        \node[state,scale=0.5](G)[left of=F]{$G$};
        \path (A)edge node{$r_{1}$}(B)
        (B)edge node{\footnotesize$r_{2}$}(C)
        (C)edge node{\footnotesize$r_{3}$}(D)
        (D)edge node{\footnotesize$r_{4}$}(E)
        (C)edge node{\footnotesize$r_{7}$}(E)
        (F)edge node{\footnotesize$r_{5}$}(A)
        (F)edge node{\footnotesize$r_{6}$}(G);
        \end{tikzpicture}
        \caption{Auxiliary graph}
         \label{Auxiliary graph}
     \end{subfigure}
     \hfill
     \begin{subfigure}{0.3\textwidth}
        \centering
        \begin{tikzpicture}[-,node distance=2cm,auto,scale = 0.25]
            \node[state,scale=0.5](A){$A$};
            \node[state,scale=0.5](B)[right of=A]{$B$};
            \node[state,scale=0.5](C)[right of=B]{$C$};
            \node[state,scale=0.5](D)[below of=C]{$D$};
            \node[state,scale=0.5](E)[below of=B]{$E$};
            \node[state,scale=0.5](F)[below of=A]{$F$};
            \node[state,scale=0.5](G)[left of=F]{$G$};
            \path (A)edge node{$r_{1}$}(B)
            (B)edge node{\footnotesize$r_{2}$}(C)
            (C)edge node{\footnotesize$r_{3}$}(D)
            (D)edge node{\footnotesize$r_{4}$}(E)
            (F)edge node{\footnotesize$r_{5}$}(A)
            (F)edge node{\footnotesize$r_{6}$}(G);
            \draw[dotted] (E) to node[swap]{\footnotesize$r_{7}$} (C);
            
            \end{tikzpicture}
            
            \caption{Updated graph}
        \label{Updated graph}
     \end{subfigure}
        \caption{Give an example to illustrate the proof idea of ~\cref{Theorem1}. $A$ to $G$ are edge demands, and $1$ to $6$ are vertices. $\{r_{1}, r_{2},..., r_{7}\}$ are multi-robotic tours. Transfer the edges of the original graph to the vertices to construct an auxiliary graph. More specifically, the vertices of the auxiliary graph are the edges of the original graph, and there exists a link in the auxiliary graph if the original two edges share the same vertices. We decrease the number of robotics by breaking the cycles in the auxiliary graph. In this example, we break the cycle in~\cref{Updated graph}, and the number of multi-robotic decreases from $7$ to $6$.}
        \label{Three graphs}
    \end{figure}

    In addition to the notation already introduced, we make the following definitions and provide the relevant properties of the SCARP.
    \begin{definition}\label{Definition1}
        The \emph{support} of a robot in a solution is the set of edges where the robot satisfies the demand
        $S_{r}(y) = \{(i,j)\in \mathcal{U},\ y^r_{ij}>0\}$
        For any solution $x,y$ we can group robots into categories based on the size of the support: 
        \begin{align*}
           & R_{0}(y) = \{ r : S_{r}(y) = \emptyset \} = \text{unused robots} \\
           & R_{1}(y) = \{ r : \lvert S_r(y)\rvert = 1 \}  = \text{robots spraying a single edge}\\
           & R_{2+}(y) = \{ r : \lvert S_r(y)\rvert \ge 2 \}  = \text{robots spraying multiple edges}
    \end{align*}
    \end{definition}
    \begin{theorem}\label{Theorem1}
        For any robotic sprayer routing problem, there exists an optimal solution $x,y$ such that  $\lvert R_{2+}(y)\rvert < m = \lvert \mathcal{U}\rvert $ that is less than $m$ robots spray multiple edges; all the remaining robots are either unused or have singleton support (spray a single edge).
    \end{theorem}
    \begin{proof}
        Assume there exists an optimal solution $x,y$ with $\lvert R_{2+}(y) \rvert \ge m$. $x$ is the set of decision variables $x_{ij}^r$, whereas $y$ is the set of actual chemical sprays $y_{ij}^r$. Transfer the edge to the vertex to construct an auxiliary graph. More specifically, the vertices of the auxiliary graph are the edges of the original graph, and there exists a link in the auxiliary graph if the original two edges share the same vertices. Then there are at least $m$ robots covering multiple edges, meaning the auxiliary graph has at least $m$ vertices and $m$ edges. Thus the auxiliary graph exists in at least one cycle. \cref{Three graphs} give an example to illustrate the proof idea.\\
        \textbf{Step 1:} There are greater than or equal to $m$ edges in the auxiliary graph, identify a cycle $C$. \\
        \textbf{Step 2:} Balancing spray around cycle $C$, update the solution by an amount sufficient to ensure one vertex in the auxiliary graph has spray 0. Note that for a cycle, every vertex visited by the cycle has exactly two adjacent edges. The divergence of the vertices equals $0$. In other words, if increase the sprays of an edge from the original graph of robot $r$, the sprays of the same edge from another robot will decline. Considering the vertices in the auxiliary graph and cycles, either of the following will happen.
        \begin{description}
        \item[\rm Case 1:] There exists an edge $(i,j)$ for which demand of spray is greater than or equal to $P$.
        %\begin{itemize} %input dot
        %    \item 
            Increasing the chemical sprays $y^r_{ij}$ = $P$, decreasing the remaining vertices sprays to $0$, and adjusting the sprays of other vertices to maintain the balance of capacity of robot $r$ and chemical demands of the vertex in the auxiliary graph. 
        %\end{itemize}
        \item[\rm Case 2:] Otherwise,
            increase the chemical sprays of a vertex to the required edge demands and decrease the same vertex sprays to $0$ for the other robots. Finally, adjust the sprays of other vertices to maintain the balance of the robot's capacity and the vertex's required demands.  
        %\end{itemize}
        \end{description}
        \textbf{Step 3:} Remove edges with zero sprays in the updated solution.\\
        \textbf{Step 4:} Cycle $C$ does not exist. If there are greater than or equal to $m$ edges in the auxiliary graph, go to Step~$1$.\\
        \textbf{Step 5:} Updating the auxiliary graph find $\bar{y}$ (updated solution of set $y$) satisfied $\lvert R_{2+}(\bar{y})\rvert < \lvert R_{2+}(y)\rvert $.\\
        
         We get a new solution $(x,\bar y)$ with fewer robots spraying multiple edges. In conclusion, there exists an optimal solution with less than $m$ robots spraying multiple edges. 
    \end{proof}
        
        Based on the above theorem, we might hope to reduce the size of the problem further by simply using robot trips that empty their whole capacity on a single row whenever that row has a demand greater than the tank capacity. The next theorem shows that this strategy does not guarantee an optimal solution.
        
        \begin{theorem}\label{Theorem2}
            For arbitrary robot spray routing problems, the condition $D_{i j} > P$ is \emph{not} sufficient to guarantee that there exists an optimal solution with a robot $r$ such that $\bar{y}^{r}_{i j}=P$.
        \end{theorem}
        \begin{proof}
            We prove this by providing a counter-example. Consider the problem in~\cref{Counterexample}.
            As shown in the figure, $D_{23}=12$ all other $D_{ij}=1$ and all $C_{ij}=1$. Let $P=8$ then $D_{23}>P$.
            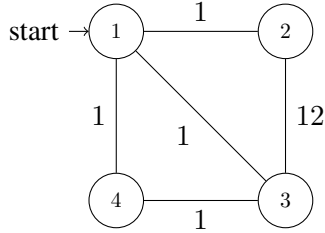
\begin{figure}[htb!]
              \centering
            \captionsetup{justification=centering,margin=1cm}
                 \begin{tikzpicture}[-,node distance=3cm,auto,scale = 0.5]        
                 \node[initial, state,scale =0.75](A){$1$};
                 \node[state,scale =0.75](B)[right of=A]{$2$};
                 \node[state,scale =0.75](C)[below of=B]{$3$};
                 \node[state,scale =0.75](D)[below of=A]{$4$};
                 \path (A)edge node{$1$}(B)
                 (B)edge node{$12$}(C)
                 (C)edge node{$1$}(D)
                 (C)edge node{$1$}(A)
                 (D)edge node{$1$}(A);
                 \end{tikzpicture}
                 \caption{Counter-example: An undirected graph and the fixed depot is vertex $1$. A robotic sprayer starts from the depot to spray chemicals to edges (1,2), (1,3), (1,4), (2,3), and (3,4). The number along each edge is the edge corresponding required demands.}
                 \label{Counterexample}
            \end{figure}
        
            \begin{table}[htb!]
            \footnotesize
                \centering
                 \caption{Optimal solution of the counter-example of ~\cref{Counterexample}}
                \begin{tabular}{c c c c c}
                    \hline
                      \multirow{2}{3em}{$r_{1}$} & $(1,2)$ & $(2,3)$ & $(3,1)$ &  \\
                     & $y^{r_{1}}_{12} = 0$ & $y^{r_{1}}_{23} =7$ & $y^{r_{1}}_{31} = 1$ & \\
                     \hline
                     \multirow{2}{3em}{$r_{2}$} & $(1,2)$ & $(2,3)$ & $(3,4)$ & $(4,1)$ \\ 
                     & $y^{r_{2}}_{12} = 1$ & $y^{r_{2}}_{23} =5$ & $y^{r_{2}}_{34} = 1$ & $y^{r_{2}}_{41} = 1$ \\
                     \hline
                \end{tabular}
                \label{Optimal solution}
                \end{table}
               \begin{table}[htb!]
                  \centering
                \footnotesize
                   \caption{New solution of the counter-example of ~\cref{Counterexample}}
                    \begin{tabular}{c c c c c c c}
                    \hline
                    \multirow{2}{1.5em}{$r_{1}$} & $(1,2)$ & $(2,3)$ & $(3,1)$&&& \\
                     & $\bar{y}^{r_{1}}_{12} = 0$ & $\bar{y}^{r_{1}}_{23} =8$ & $\bar{y}^{r_{1}}_{31} = 0$ &&& \\
                     \hline
                     \multirow{2}{1.5em}{$r_{2}$} &$(1,2)$ & $(2,3)$ & $(3,4)$ & $(4,1)$ & $(1,3)$ &$(3,1)$\\
                     & $\bar{y}^{r_{2}}_{12} = 1$ & $\bar{y}^{r_{2}}_{23} =4$ & $\bar{y}^{r_{2}}_{34} = 1$ & $\bar{y}^{r_{2}}_{41} = 1 $ & $\bar{y}^{r_{2}}_{13} = 1$ & $\bar{y}^{r_{2}}_{31} = 0$\\
                     \hline
                \end{tabular}
                \label{New solution}
           \end{table}
          %\newpage
            \cref{Optimal solution} gives an optimal solution value to SCARP formulation 1 of $7$. \cref{New solution} shows when the first robot is fixed to only spray edge $(2,3)$, the optimal solution cost will be increased to $9$. Thus, for any robot sprayer routing problem, if $D_{i j} > P$ for some $(i,j)\in \mathcal{U}$ and here not always existing an optimal solution $\bar{x},\bar{y}$ and a robot $r$ such that $\bar{y}^{r}_{i j}=P$. 
        \end{proof}
        
\subsection{SCARP formulation 2: large edge demands formulation}
        \label{SCARP}
        Based on \protect\cref{Theorem1} there is an opportunity to create a simpler formulation when some of the edges' demands $D_{ij}$ are much greater than $P$. 
        
        Apply \emph{Dijkstra's algorithm} to calculate the shortest path distance from the start vertex $1$ (depot) to other vertices in the graph. Let $SP_{1 i}$ be the shortest path distance from vertex $1$ to vertex $i$. $\mathcal{R}_{new}=\{ 1,2,...,$min$(\lvert \mathcal{U}\rvert -1, \lvert \mathcal{R}\rvert  = \frac{\sum_{(i,j)\in \mathcal{U}}D_{i j}}{P}) \}$ set of robots that work on multiple edges. $z_{i j} \in  \mathbb Z_{+}$ the number of singleton robots that spray edge $(i,j)\in \mathcal{U}$.
%   \newpage
   \begin{mini!}|s|[1]
            {}{\sum_{(i,j)\in\mathcal{U}}z_{i j}(SP_{1 i} + C_{i j} + SP_{j 1}) +\sum_{r\in \mathcal{R}_{new} }\sum_{(i,j)\in \mathcal{U}}C_{i j} (x_{i j}^r+x_{j i}^r)  \quad (j,i)\in\mathcal{E} }{}{\label{Const2-0}}
            \addConstraint{\sum_{(i,j)\in\mathcal{U}} y_{i j}^r}{\leq P}{\quad \forall r\in \mathcal{R}_{new}\label{Const2-1}}
            \addConstraint{y_{i j}^r}{\geq 0}{\quad \forall(i,j)\in\mathcal{U}, r\in\mathcal{R}_{new}\label{Const2-2}}
            \addConstraint{\sum_{r \in \mathcal{R}}y_{i j}^r + z_{i j}P}{\ge D_{i j}}{\quad\forall(i,j)\in \mathcal{U}\label{Const2-3}}
            \addConstraint{y_{i j}^r}{\leq P(x_{i j}^r+x_{j i}^r)}{\quad\forall(i,j)\in \mathcal{U},r \in \mathcal{R}_{new}\label{Const2-4}}
            \addConstraint{\sum_{i\in\mathcal{I}_k}x_{i k}^r}{=\sum_{j\in \mathcal{Q}_k}x_{k j}^r}{\quad\forall k\in \mathcal{V}, r\in \mathcal{R}_{new}\label{Const2-5}}
            \addConstraint{\sum_{i\in \mathcal{S},j\in\mathcal{T}}Px_{i j}^r}{\geq \sum_{i,j\in \mathcal{T}}y_{i j}^r}{\quad \substack{\forall \mathcal{S}\subset \mathcal{V},S\cup\mathcal{T}=\mathcal{V},\\ 1\in\mathcal{S}\mathcal{T}=\mathcal{V}\setminus{\mathcal{S}}, r \in\mathcal{R}_{new}}\label{Const2-6}}
            \addConstraint{x_{i j}^r}{\in\{0,1\}}{\quad\forall(i,j)\in \mathcal{E},r \in \mathcal{R}_{new}\label{Const2-7}}
            \addConstraint{z_{i j}}{\le \frac{D_{i j}}{P}}{\quad\forall(i,j)\in \mathcal{U}\label{Const2-8}}
       \end{mini!}

       The key idea is splitting the large demand on a single edge amongst vehicles to involve fewer vehicles to decrease the number of decision variables $x^r_{i j}$ and $y^r_{i j}$ when $\lvert\mathcal{R}\rvert > \lvert\mathcal{U}\rvert-1$. \cref{Theorem1} indicates that there exits at least one robotic sprayer robot $r$ that $\bar{y}^{r}_{i j}=P$ when $\lvert\mathcal{R}\rvert > \lvert\mathcal{U}\rvert-1$. More specifically, there are at most $\lvert\mathcal{U}\rvert-1$ multi-robotic and at least  $\lvert\mathcal{R}\rvert- ( \lvert\mathcal{U}\rvert-1)$ single-robotic (Except in the case of the last robot covering the rest of the spray). The objective function ~(\ref{Const2-0}) consists of singleton and multi-robotic' costs. Similarly, Constraints~ (\ref{Const2-3}) indicate each edge can be serviced by singleton and multi-robotic sprayers. Hence, the total edge demand comes from the chemicals sprayed by singleton and multi-robotic sprayers.

\section{Methodology}\label{sec5}

\subsection{Lazy-constraints}
        \noindent Lazy constraints are useful when the full set of constraints is too large to include in the initial formulation. The main idea is to remove the costly constraints and get a relaxed problem with the remaining constraints. And then, we reintroduce that constraint if the relaxation solutions violate some of the constraints we left out.
        
        The basic formulation described in \cref{sec3} 
        is costly for a large-scale problem. Because the number of Constraints~\eqref{Const6}  grows exponentially as the vertices increase in the graph. We simply cannot afford to list them all for decent size graphs. An alternative is to treat them as lazy constraints and only add them when they are violated. 
        
        The main steps for applying lazy constraints are as follows:
      \begin{itemize}
            \item Remove the connectivity  Constraints \eqref{Const6};
            \item For any LP solution obtained during the branch-and-bound, test if the solution violates Constraints \eqref{Const6} for any robot $r$. If it does, go to the next step, otherwise, skip the lazy constraints;
            \item Use the \emph{BFS} algorithm to generate suitable subset $\mathcal{S}$ of $\mathcal{N}$. Pass the $\mathcal{S}$ to build cuts in the branch-and-bound process.
           
        \end{itemize}
        
       The constraints could be added at both branch-and-bound nodes with integer and fractional LP solutions. For the sake of computational efficiency, we add constraints for fractional solutions early in the search and all infeasible integer solutions. More specifically, we add lazy constraints on a limited number of nodes (not too deep in the branch and bound tree) when the gap between the lower and upper bound is not too small.

\subsection{Symmetry elimination constraints} 
        There are many equivalent ways to encode the same  solution in our formulation simply by renumbering the robots. We can eliminate this symmetry with Constraints~\eqref{3a}. The Constraints~\eqref{3a} order the robots from less costly to more to eliminate symmetric solutions. 
        \begin{equation}
            \tag{3a}
            \sum_{(i,j) \in \mathcal{A}} C_{ij}x_{ij}^r \ge \sum_{(i,j) \in \mathcal{A}} C_{ij}x_{ij}^{r-1} \hspace{1em} \forall r\in\mathcal R_{new}:r>1 \label{3a}
        \end{equation}
        
        Another form of symmetry is due to the fact that traversing a route in the reverse order does not change the solution. To eliminate this symmetry, let $N^i_{1} = \{i \in (i,1): \forall (i,1) \in \mathcal{I}_1\}$ be the set of vertices for arcs coming into start vertex $1$. The Constraints~(\ref{3b}) orients the route by not allowing the return arc at the depot to connect to a lower indexed node than the outgoing arc. However, this doesn't help if the first and last edge visited by the robot are the same.
        % No paragraph break here
        \begin{equation}
            \tag{3b}
            \sum_{i = 1}^{k} x_{1i}^r \ge x_{k1}^r \hspace{1em} \forall k \in  N^i_{1} , \hspace{0.5em}  \forall r\in\mathcal R_{new} \label{3b}
        \end{equation}
        
\subsection{Heuristic repair method}
        \noindent We provide a heuristic repair method that can be applied to any infeasible (fractional) solutions. The main idea is to construct heuristic solutions based on infeasible solutions and use them as new incumbents during the branch-and-bound.
        
        We detail its main components in the following paragraphs.
        
        \textbf{Preprocessing step:} Compute shortest path and corresponding shortest \emph{distance} (cost) from vertex $i$ to vertex $j$ ($\forall\ i,j \in \mathcal{V}$) by \emph{Dijkstra’s algorithm}, represented by $SPath_{ij}$ and $SP_{ij}$, respectively.
      
        Using an infeasible solution with feasible $y$ values to construct a connected path that covers all edges $(i,j)$ such that $y_{ij}^r > 0$ by the \emph{Greedy routing algorithm} ($y_{ij}^r$ can be fractions, thus they are always feasible). It is worth noting that in chaining paths from one spray edge to the next we may be traversing an arc in the same direction twice and that we can trivially improve the solution by removing such duplicates. Then generate $x_{ij}^r$ according to the connected path. At last, sort $x_{ij}^r$ and $y_{ij}^r$ according to the symmetry elimination constraints in \eqref{3a} and \eqref{3b}.
        \RestyleAlgo{ruled}
        \SetAlFnt{\small}
        \SetKwComment{Comment}{$\triangleright$\ }{ $\triangleright$\ }
        \begin{algorithm}[htp]
        \caption{Greedy routing algorithm} \label{Alg.1}
        \KwData{Infeasible solutions with feasible $y = \{y_{ij}^r : (i,j) \in \mathcal{U}, r \in \mathcal{R}\}$, $SPath_{ij}$ , $SP_{ij}$}
        \KwResult{ Feasible solutions $x = \{ x_{ij}^r : (i,j) \in \mathcal{E}, r \in \mathcal{R}\}$, $y$\;}
        $F\_edge = \{(i,j) : y_{ij}^r > 0\}$\;
        $B\_edge = \{(j,i) : y_{ij}^r > 0\}$\;
        $All\_edge = append!(F\_edge,B\_edge)$\;
        $s = 1$ \Comment*[r]{Depot}
        \SetKwProg{Fn}{Function}{:}{end}
        
         \SetKwFunction{FMain}{Greedy routing}
        \Fn{\FMain{$y$, $SPath_{ij}$, $SP_{ij}$}}{
        $route =$ Connectable path $(y, SPath_{ij}, SP_{ij})$\;
        $path =$ Remove duplicate $(y, SPath_{ij}, SP_{ij}, route)$\;
        \For{$(i,j) \in path$}{
        $ x \gets x_{ij}^r = 1$\;}
        $order  \gets$ \ref{3a}, \ref{3b}\;
        Sorting $x , y \gets order$\;
        \Return $x , y$
         }

        \SetKwFunction{FMain}{Connectable path}
        \Fn{\FMain{$y$, $SPath_{ij}$, $SP_{ij}$}}{
        $ route \gets  \emptyset$ \;
        \While{$\lvert All\_edge \rvert \neq 0$}{
            $n \gets $ vertices passed by $F\_edge$ \;
            $v \gets$ min($SP_{sn}$ , n)\;
          \If{$\lvert SPath_{sv}\rvert  \ge 2 $}{
           $ route \gets$ append $SPath_{sv}$\;}
           $ pass\_vertex \gets \{j: (v,j) \in F\_edge\} $\;
          \If{ $\lvert \text{pass\_vertex}\rvert  = 0 $}{
           $ reverse\_vertex \gets \{j: (j,v) \in B\_edge\} $\; 
           $ pass\_vertex \gets$ append $reverse\_vertex[1]$\;
            }
            $ route \gets $append $pass\_vertex[1]$\;
            update $F\_edge$, $B\_edge$  by removing covered edges in both directions\;
            $All\_edge \gets append! ( F\_edge, B\_edge )$\;
            }
           \If{route[end] != {1}}{
           $ back = SPath_{route[end]s}$ \Comment*[r]{ find shortest path from $route[end]$ to $s$ }
           $ route \gets $append $SPath_{route[end]s}$\;
           } 
           \Return $route$
        }
        \SetKwFunction{FMain}{Remove duplicate}
        \Fn{\FMain{$y$, $SPath_{ij}$, $SP_{ij}$, $route$}}{
        \For{$i \gets 1$  \KwTo $\lvert route\rvert - 1$ }{
        $ path \gets {(route[i], route[i+1])}$\;}
        \If{path has duplicated edges}{ 
        remove duplicated edges from $path$
        }
        \Return path
        }
       
        \end{algorithm}

        In Alg.~\ref{Alg.1} Greedy routing function, some edges need to be reversed, such as reversing $(i,j)$ to $(j,i)$ to get a connectable path when the outgoing vertex is $j$. The traversal of $(i,j)$ or $(j,i)$ are both seen as a coverage of $(i,j)$. If it is an edge in $F\_edge$,  remove it from  $F\_edge$, then swap the edge order and remove it from $B\_edge$, and vice versa.   
        
        The heuristic repair method also generates many rejected solutions, which return worse objective values than the exits. To improve the efficiency, we set a parameter $\Gamma$ representing the number of times the heuristic will be run without improvement before abandoning any further attempts to run the repair heuristic. Section~\ref{sec7} provides an analysis of the parameter $\Gamma$.
    
\section{Numerical results}\label{sec6}

    In this section, we test the effectiveness of our algorithms for solving this agricultural vehicle spraying problem. We consider three main classes of instances. The first and second classes of instances are generated by ourselves based on the attributes of the real-world examples described by \cite{plessen2019optimal}.

    \begin{table*}[tb]
    \centering
    \caption{Comparative results of the Basic and the Large edge demands models on the same graph LD with different edge demands. LD consists of 16 vertices and 22 edges, which is displayed in~\cref{Fig.4}. Robot capacity equals $10$.}
    \label{Computational results between Basic model and High models}
    \begin{adjustbox}{max width=\textwidth}
    \begin{tabular}{l l l l l l l l l l}
    \hline
    Instance  & sum(D) & B\_r & L\_r & B\_obj & L\_obj & B\_gap & L\_gap & B\_t & L\_t \\ \hline
    LD\_1  & 260.131 & 27 & 5 & 6110.93  & 6110.93  & 37.07\% & 0.96\% & 7200 & 7006.25 \\ 
    LD\_2 & 390.195 & 40 & 5 & 13675.18 & 13572.47 & 40.75\% & 0.51\% & 7200 & 6358.96 \\ 
    LD\_3  & 520.26  & 53 & 4 & 24058.06 & 23991.43 & 40.02\% & 0.66\% & 7200 & 6014.61 \\ \hline
    \end{tabular}
    \end{adjustbox}
    \end{table*}
    
    The first class of large edge demands instances LD\_1, LD\_2, and LD\_3 contains 16 vertices and 22 edges. The first set of instances comprise an orchard area with different spraying requirements. These instances are applied to make a comparison between the basic and large edge demands formulations. The second class of instances contain seven different graphs with or without random obstacles. Instances A, B, and C are small-scale orchards without obstacles, whereas instances E, F, and G have one or two random obstacles. The largest size graph G has 98 vertices and 245 edges. We also test our models and algorithms on a set of classical CARP instances, which contains 34 undirected graphs designed by \cite{belenguer2003cutting}. These instances are defined on 10 different graphs by changing capacity of the vehicles. In classical CARP instances set, all edges are required. The algorithm described in this paper was coded in Julia 1.7 using Gurobi 9.5 as a MIP solver, and tested on computers with 2.70~GHz Intel Xeon-Gold-6150 processors. All computation times reported in this paper are in seconds with a limit of 7200 seconds unless otherwise specified. 
    
     \begin{figure}[htb!]
        \centering
         \includegraphics[scale=0.5]{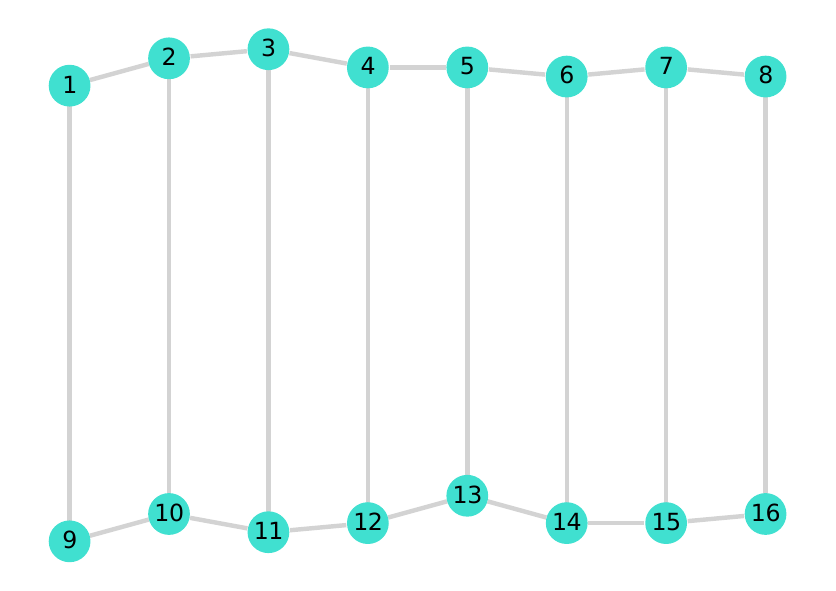}
         \caption{ The graph of LD: The instance LD has 16 vertices and 22 edges. We are changing the edge demand along each edge to get LD\_1, LD\_2, and LD\_3 in~\cref{Computational results between Basic model and High models}. }
         \label{Fig.3}
         \end{figure}
    
    \cref{Computational results between Basic model and High models} shows the results for the first class of large demands instances (The detailed results are provided in~\cref{ Details of computational results between Basic model and High edge-demand model} at Appendix~\ref{secA1} ). The column ``sum(D)'' gives the total spraying demand for each instance. The column ``B\_r'' is the number of required robots of the basic robotic spraying formulation. The column ``L\_r'' shows the number of multi-robotics (robots spraying multiple edges, mentioned in \cref{Definition1}) in the large edge demands formulation. ``B\_obj'' and ``L\_obj'' report the objective solution found by the basic and large edge demands model. The ``B\_gap'' and ``L\_gap'' are radio percentages with respect to the upper and lower bound of the basic and large edge demands models. ``B\_t'' and ``L\_t'' give the computing time of the two models (running times in seconds). 
    
\subsection{Computational results for Large edge demands instances}
    \noindent  LD\_1, LD\_2, and LD\_3 are instances with different edge demands.  There exist some edges with requirements greater than the robot's capacity.  From \cref{Theorem1}, we know the large edge demands model works well if there have singleton robotic sprayers. \cref{Computational results between Basic model and High models} reports the computational results of the basic and large edge demands model with respect to these instances. The algorithms are implemented with 8 threads per CPU task. The Gurobi parameter random number ``seed'' acts as a small perturbation to the solver, and typically leads to different solution paths. We set the seed  parameter to the values 1, 2, and 3. We displayed the average results from the three rounds.
    
    \cref{Computational results between Basic model and High models} shows the large edge demands model outperforms the basic model regarding solution quality and has a faster computational time. A lower computational time for the large edge demands model indicates that there has fewer robot routes to optimize, which reduces the total number of variables. Besides, the basic model returns a big gap, whereas the  large edge demands model reports a relatively small gap in all three instances. It is worth noting that parameter ``seed'' numbers lead to different solution paths, which may affect the average calculation. More specifically, the large edge demands reaches an optimal solution before 7200 seconds in some seed numbers, but some seed paths can not find the optimal solution in 7200 seconds.  Thus the time is less than 7200 seconds and returns a small gap for the  large edge demands model in~\cref{Computational results between Basic model and High models}. 
    
\begin{figure}[htb]
    \centering
     \includegraphics[scale=0.4]{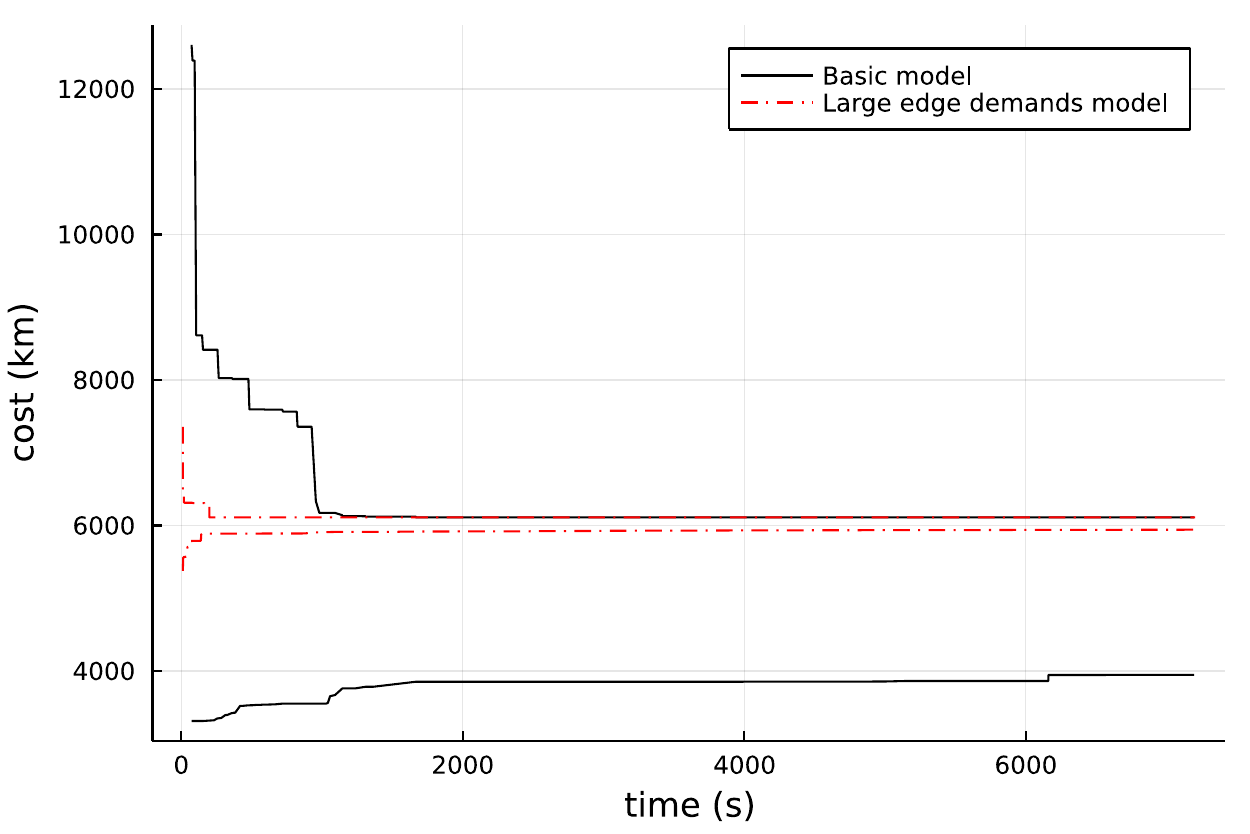}
     \caption{ LD\_1 comparison results: The black solid  and red dash-dotted lines are the upper and lower bounds of the basic and large edge demands model, respectively. We set seed =1, time limit =7200 seconds. The graph shows that the large edge demands model (red dash-dotted lines) has a fast convergence speed and returns a small gap.}
     \label{Fig.4}
     \end{figure}

\begin{figure}[htb]
    \centering
     \includegraphics[scale=0.4]{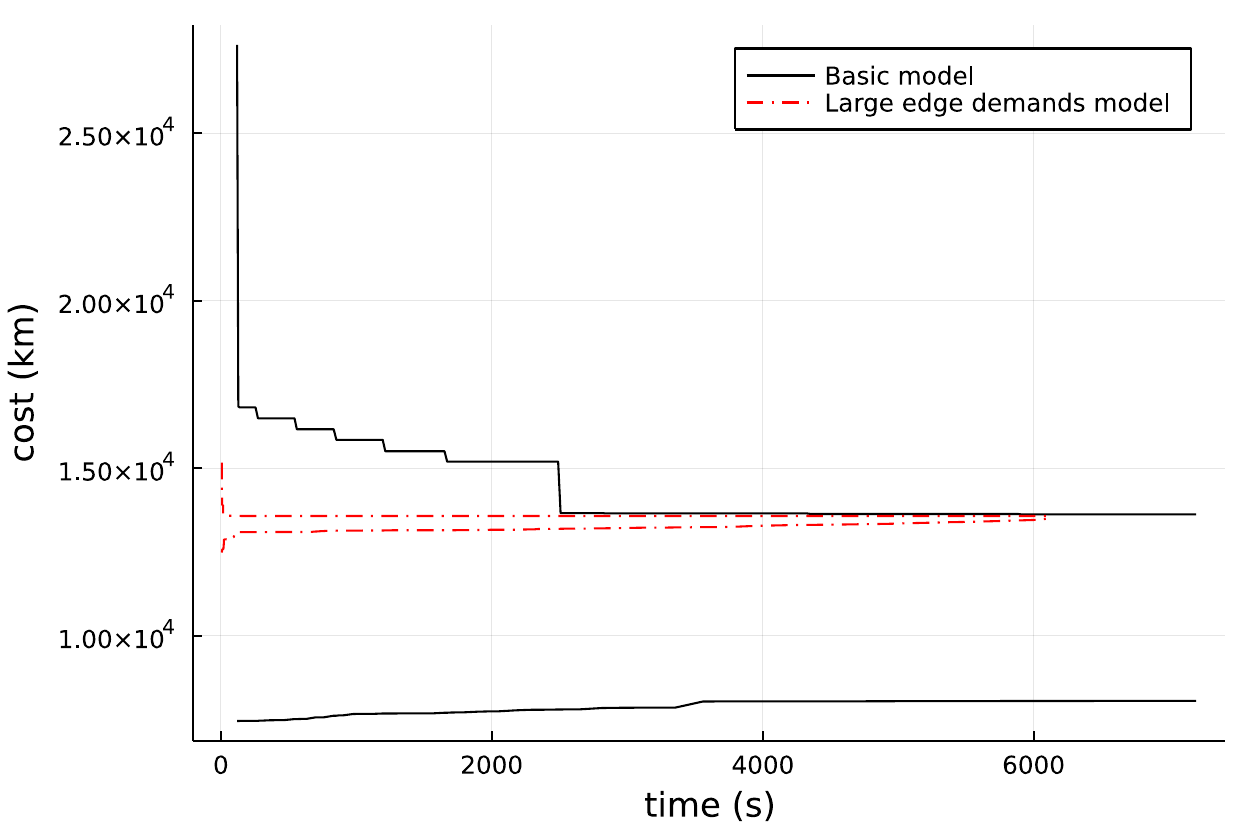}
     \caption{ LD\_2 comparison results: The black solid and red dash-dotted lines are the upper and lower bounds of the basic and large edge demands model, respectively. The graph shows that the large edge demands model (red dash-dotted lines) converges quickly and reaches optimal at about 6000 seconds, while the basic model still returns a big gap.}
     \label{Fig.5}
     \end{figure}

\begin{figure}[ht!]
        \centering
         \includegraphics[scale=0.4]{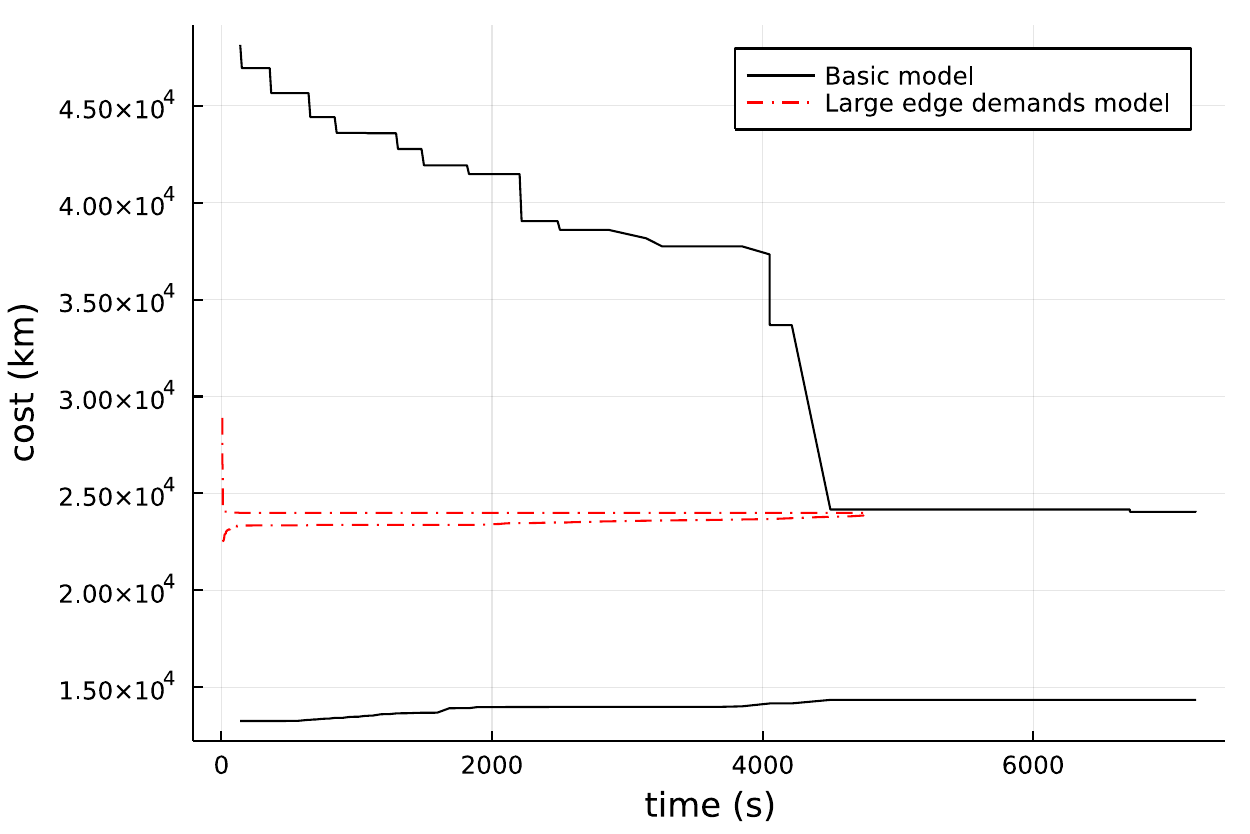}
         \caption{ LD\_3 comparison results: The black solid and red dash-dotted lines are the upper and lower bounds of the basic and large edge demands model, respectively. The graph shows that the large edge demands model (red dash-dotted lines) has a better root solution, whereas the basic model has a significant gap in the root node. The large edge demands model finds an optimal solution at about 4500 seconds, while the basic model returns a big gap.}
         \label{Fig.6}
         \end{figure}

    Figures \ref{Fig.4}, \ref{Fig.5} and \ref{Fig.6} show the large demand instances' performances in the branch-and-bound tree. The $x$ axis is time (s), and the $y$ axis is the value of the lower and upper bounds. The black solid line goes to the basic model, while the red dash-dotted lines represent the large edge demands model. The upper black solid and red dash-dotted lines show the changes in the upper bound of the basic and large edge demands model, respectively. Similarly, the bottom black and red dash-dotted lines reflect the changing of lower bounds. We use a consistent random number seed = 1.
    
\subsection{Computational results for different techniques}
      
    The algorithms in \cref{Computational results of different techniques} are described in \cref{sec3} and \cref{sec5}. The first model is the basic model, denoted as the ``Basic-model''. We add the lazy constraints to the basic model to get the result of ``Lazy-constraints''. The symmetry elimination constraints continue to be added to obtain the result of ``Sym-elimination''. Finally, we add the heuristic repair method to get the result of the``Heuristic-repair''. Besides, we check the performance of symmetry eliminations by getting rid of symmetry constraints from the heuristic repair method (``Heuristic-noSym''). The experiments are implemented with 32 threads per CPU task. We display the average data from the three rounds of different ``seed'' numbers. Different seed numbers generally lead to different solution paths and return similar solutions. However, it sometimes returns significant differences in time and solutions, remarkably affecting average values. The ``obj'' returns the objective solution. At the same time, the column ``gap'' is the deviation between the upper and lower bound. The column ``t'' displays the computing time in seconds, and the parameter ``$\lvert Acc \rvert$'' represents the number of accepted heuristic repair solutions.
    
    \cref{Computational results of different techniques} shows that the basic model only works for small-scale problems A and B. In contrast, lazy constraints and symmetry eliminations could solve medium-sized problems E and F. The heuristic repair method shows excellent solution quality and can be applied to large instances G. Computational time indicate that the Basic model is the most expensive one. At the same time, the other algorithms look similar to time consumption. We tested all of the algorithms in 7200 seconds (2 hours). After running 20 hours of instance G, we get a similar solution to the 2 hours solution, which shows the longer running time doesn't improve the solution a lot.  
     
    \begin{sidewaystable*}[p]
        \caption{Comparative results of different techniques on the seven undirected instances from $A$ to $G$. This set was designed based on the real-world orchard data \citep{plessen2019optimal}. The column ``sum(D)'' is the total edge demands for the corresponding instance. Columns ``obj'' and ``gap'' show the optimal cost and deviation ($gap(\%)  = 100 (UB-LB) / LB$, where $UB$ is the solution cost, and $LB$ is the lower bound ). The column ``$\lvert Acc \rvert$'' indicates the number of accepted heuristic solutions. The best performances are highlighted in grey cells, whereas the bold are the worst solution for each instance. \cref{Details of computational results of different techniques} in Appendix~\ref{secA1} provides detailed results.}
        \label{Computational results of different techniques}
        % \footnotesize
        \begin{adjustbox}{max width=\textwidth}
        \begin{tabular}{l l l l l l l l l l l l l l l l l l l l l l l l  }
        \hline
        \multirow{2}{*}{Ins.} &
          \multirow{2}{*}{$\lvert N \rvert$} &
          \multirow{2}{*}{$\lvert E \rvert$} &
          \multirow{2}{*}{p} &
          \multirow{2}{*}{r} &
          \multirow{2}{*}{sum(D)} &
          \multicolumn{3}{c}{Basic-model} &
          \multicolumn{3}{c}{Lazy-constraints} &
          \multicolumn{3}{c}{Sym-elimination} &
          \multicolumn{4}{c}{Heuristic-repair} &
          \multicolumn{4}{c}{Heuristic-noSym} \\ \cline{7-23} 
         &
           &
           &
           &
           &
           &
          \multicolumn{1}{c}{obj} &
          \multicolumn{1}{c}{gap} &
          \multicolumn{1}{c}{t (s)} &
          \multicolumn{1}{c}{obj} &
          \multicolumn{1}{c}{gap} &
          \multicolumn{1}{c}{t (s)} &
          \multicolumn{1}{c}{obj} &
          \multicolumn{1}{c}{gap} &
          \multicolumn{1}{c}{t (s)} &
          \multicolumn{1}{c}{obj} &
          \multicolumn{1}{c}{gap} &
          \multicolumn{1}{c}{t (s)} &
          \multicolumn{1}{c}{$\lvert Acc \rvert$} &
          \multicolumn{1}{c}{obj} &
          \multicolumn{1}{c}{gap} &
          \multicolumn{1}{c}{t (s)} &
          \multicolumn{1}{c}{$\lvert Acc \rvert$} \\ \hline
        \begin{tabular}[c]{@{}l@{}}Ins.A\end{tabular} &
          16 &
          22 &
          20 &
          1 &
          13.01 &
          32.12 &
          0.00\% &
          \colorbox{gray}{0.85} &
          32.12 &
          0.00\% &
          2.11 &
          32.12 &
          0.00\% &
          1.98 &
          32.12 &
          0.00\% &
          \textbf{5.69} &
          1 &
          32.12 &
          0.00\% &
          4.26 &
          1 \\ 
        \begin{tabular}[c]{@{}l@{}}Ins.B\end{tabular} &
          20 &
          28 &
          20 &
          1 &
          17.14 &
          42.64 &
          0.00\% &
          \textbf{28.01} &
          42.64 &
          0.00\% &
          2.07 &
          42.64 &
          0.00\% &
          \colorbox{gray}{2.04} &
          42.64 &
          0.00\% &
          5.97 &
          1 &
          42.64 &
          0.00\% &
          4.03 &
          1 \\ 
        \begin{tabular}[c]{@{}l@{}}Ins.C\end{tabular} &
          30 &
          43 &
          20 &
          2 &
          28.36 &
          fail &
           &
           &
          76.15 &
          0.00\% &
          \colorbox{gray}{49.99} &
          76.15 &
          0.00\% &
          62.07 &
          76.15 &
          0.00\% &
          79.04 &
          3 &
          76.15 &
          0.00\% &
          \textbf{91.90}&
          3 \\ 
        \begin{tabular}[c]{@{}l@{}}Ins.D\end{tabular} &
          42 &
          61 &
          20 &
          2 &
          32.58 &
          fail &
           &
           &
          84.30&
          0.00\% &
          5.74 &
          84.30 &
          0.00\% &
          \colorbox{gray}{5.62} &
          84.30 &
          0.00\% &
          \textbf{10.93} &
          6 &
          84.30 &
          0.00\% &
          9.23 &
          7 \\ 
        \begin{tabular}[c]{@{}l@{}}Ins.E\end{tabular} &
          56 &
          82 &
          20 &
          3 &
          46.59 &
          fail &
           &
           &
          127.65 &
          2.68\% &
          7200 &
          127.65 &
          \textbf{3.62\%} &
          7200 &
          127.65 &
          \colorbox{gray}{2.45\%} &
          7200 &
          8 &
          127.65 &
          2.75\% &
          7200 &
          6 \\ 
        \begin{tabular}[c]{@{}l@{}}Ins.F\end{tabular} &
          78 &
          115 &
          20 &
          3 &
          52.98 &
          fail &
           &
           &
          150.74 &
          6.91\% &
          7200 &
          151.98 &
          7.42\% &
          7200 &
          152.19 &
          \textbf{7.64\%} &
          7200 &
          8 &
          149.79 &
          \colorbox{gray}{6.06\%} &
          7200 &
          8 \\ 
        \begin{tabular}[c]{@{}l@{}}Ins.G\end{tabular} &
          98 &
          145 &
          20 &
          4 &
          74.78 &
          fail &
           &
           &
          fail &
           &
           &
          fail &
           &
           &
          235.77 &
          \colorbox{gray}{15.75\%} &
          7200 &
          13 &
          251.58 &
          \textbf{20.72}\% &
          7200 &
          6 \\ \hline
        \end{tabular}
        \end{adjustbox}
    \end{sidewaystable*}

    The results of small-scale instances A and B show that all the algorithms could find an optimal solution, but it is considerably costlier for the heuristic repair method. This is because the heuristic repair method spends much time generating feasible solutions. As the size of the instances increases, the basic model fails to find an optimal solution in instances C and D, while other techniques work well. Instances E, F, and G show that the heuristic repair method works better in large-scale problems. Except for the heuristic technique, all the other methods failed to get solutions for the large-scale instance G. \cref{Computational results of different techniques} shows that symmetry eliminations do not always work well in our instances. It saves time in test C and improves the solution quality for test E and G, but not helpful for the others. Our future work is to make symmetry eliminations work stable for our algorithms. As mentioned in~\cref{sec5}, the heuristic repair method involves the parameter $\Gamma$, which affects the method's overall performance. We reported in~\cref{sec4}. We set $\Gamma$ = 1000 in this set of tests. \cref{Details of computational results of different techniques} in Appendix~\ref{secA1} provides detailed results. 

\subsection{Computational results for benchmark instances}
        \begin{table*}[tbh!]
        \centering
        \caption{Comparative results of the heuristic repair method with state-of-art algorithms from ~\cite{chen2016hybrid} on the 34 undirected instances of \emph{val} set. This set was first designed to find CARP's lower bound by~\cite{benavent1992capacitated}. The column ``Best-Known'' is the best solution from previous literature. ``LB'' and ``LB-H'' show the lower bounds from the literature and this paper's heuristic repair method, respectively. The better solutions are highlighted in grey cells, whereas the star indicates the solution matches the ``Best-Known'' solution.}
        \label{Computational results of "val" file}
        % \footnotesize
        \begin{adjustbox}{max width=0.75\textwidth}
        \begin{tabular}{llllllllll}
        \hline
        File   & $\lvert N \rvert$ & $\lvert E \rvert$ & r  & p   & sum(D) & LB  & LB-H  & Best-Known   & Heuristic-repair  \\ \hline
        val1a  & 24  & 39  & 2  & 200 & 358    & 173 & 173   & 173 & 173*                  \\ 
        val1b  & 24  & 39  & 3  & 120 & 358    & 173 & 173   & 173 & 173*                \\ 
        val1c  & 24  & 39  & 8  & 45  & 358    & 235 & 191   & 245 & \colorbox{gray}{237}                  \\ \hline
        val2a  & 24  & 34  & 2  & 180 & 310    & 227 & 227   & 227 & 227*                  \\ 
        val2b  & 24  & 34  & 3  & 120 & 310    & 259 & 259   & 259 & 259*                \\ 
        val2c  & 24  & 34  & 8  & 40  & 310    & 455 & 417   & 457 & \colorbox{gray}{455}               \\ \hline
        val3a  & 24  & 35  & 2  & 80  & 137    & 81  & 81     & 81      & 81*          \\ 
        val3b  & 24  & 35  & 3  & 50  & 137    & 87  & 87     & 87      & 87*          \\ 
        val3c  & 24  & 35  & 7  & 20  & 137    & 137 & 115   & 138     & 138*        \\ \hline
        val4a  & 41  & 69  & 3  & 225 & 627    & 400 & 400   & 400     & 400*        \\ 
        val4b  & 41  & 69  & 4  & 170 & 627    & 412 & 382   & 412     & 418      \\ 
        val4c  & 41  & 69  & 5  & 130 & 627    & 428 & 382   & 428     & 434     \\ 
        val4d  & 41  & 69  & 9  & 75  & 627    & 520 & 437   & 528     & 536      \\ \hline
        val5a  & 34  & 65  & 3  & 220 & 614    & 423 & 414   & 423     & 423*       \\ 
        val5b  & 34  & 65  & 4  & 165 & 614    & 446 & 428   & 446     & 446*        \\ 
        val5c  & 34  & 65  & 5  & 130 & 614    & 469 & 433   & 474     & 474*   \\ 
        val5d  & 34  & 65  & 9  & 75  & 614    & 571 & 499   & 575    & 579   \\ \hline
        val6a  & 31  & 50  & 3  & 170 & 451    & 223 & 223   & 223     & 223*       \\ 
        val6b  & 31  & 50  & 4  & 120 & 451    & 231 & 222   & 233     & 233*        \\ 
        val6c  & 31  & 50  & 10 & 50  & 451    & 311 & 239   & 317     & 317*     \\ \hline
        val7a  & 40  & 66  & 3  & 200 & 559    & 279 & 279   & 279     & 279*     \\ 
        val7b  & 40  & 66  & 4  & 150 & 559    & 283 & 283   & 283     & 283*       \\ 
        val7c  & 40  & 66  & 9  & 65  & 559    & 333 & 272   & 334     & 334*  \\ \hline
        val8a  & 30  & 63  & 3  & 200 & 566    & 386 & 386   & 386     & 386*       \\ 
        val8b  & 30  & 63  & 4  & 150 & 566    & 395 & 395   & 395     & 395*        \\ 
        val8c  & 30  & 63  & 9  & 65  & 566    & 517 & 439   & 521     & 527     \\ \hline
        val9a  & 50  & 92  & 3  & 235 & 654    & 323 & 320   & 323     & 323*      \\ 
        val9b  & 50  & 92  & 4  & 175 & 654    & 326 & 322   & 326     & 327       \\ 
        val9c  & 50  & 92  & 5  & 140 & 654    & 332 & 289   & 332     & 357      \\ 
        val9d  & 50  & 92  & 10 & 70  & 654    & 382 & 327   & 389     & 578      \\ \hline
        val10a & 50  & 97  & 3  & 250 & 704    & 428 & 425   & 428     & 428*       \\ 
        val10b & 50  & 97  & 4  & 190 & 704    & 436 & 432   & 436     & 436*    \\ 
        val10c & 50  & 97  & 5  & 150 & 704    & 446 & 428   & 446     & 527          \\ 
        val10d & 50  & 97  & 10 & 75  & 704    & 524 & 459   & 525     & 649              \\ \hline
        \end{tabular}
         \end{adjustbox}
        \end{table*}

        \noindent The set (\emph{val} benchmark data set) contains 34 undirected instances designed by \cite{benavent1992capacitated} to find the lower bounds of the CARP. These instances have been generated by varying the capacity of the vehicles for each graph (denoted val1a,valb,...,val10d). \cref{Computational results of "val" file} reports the results for \emph{val} benchmark data set. The columns ``$\lvert N \rvert$'' and ``$\lvert E \rvert$'' are the number of vertices and edges. The ``r'' is the number of vehicles, and ``p'' is the capacity of the identity vehicles. The ``LB'' is the lower bound are obtained by \cite{belenguer2003cutting}. The column ``Best-Known'' are the best results reported in previous CARP literature (summarized by \cite{chen2016hybrid}). \cite{chen2016hybrid} presented a hybrid metaheuristic approach to solve the CARP. It is worth noting that all the results from other literature are not allowed to spit the demand along each required edge. The columns ``LB-H'' and ``Heuristic-repair'' are lower bounds and objective solutions computed by the heuristic repair method. The experiments are implemented with 32 threads per CPU task, seed = 1, $\Gamma$ = 4000.
        
        \cref{Computational results of "val" file} show that the heuristic repair method finds a better solution in instances ``val1c'' and ``val2c'' since our problem allows splitting the demand on a single demand edge amongst robots. Some latest works with these benchmark sets did not provide a detailed solution or generated new data sets based on \cite{benavent1992capacitated} benchmark sets for an extension CARP. \cite{arakaki2019efficiency} presented an efficiency-based path-scanning heuristic for the CARP, which introduced the concept of the efficiency rule. Arakki and Usberti's results outperformed the other algorithms in average time but return a worse deviation (the average deviation comes from the lower bounds $GAP(\%)  = 100 (UB-LB) / LB$, where $UB$ is the solution cost, and $LB$ is the lower bound) of \emph{val} instances set. This paper focuses on the comparison of optimal cost results for \emph{val} set. Thus, we make a comparison with \cite{chen2016hybrid}. We set the time limit to 2 hours and got an improved solution in \emph{val4B} from 418 to 412 after extending the time limit to 3 hours. However, after increasing the time limit to 3 hours, the lower bound of most instances is slightly raised and returns the same objective solutions as 2 hours.   
       
\section{Parameter analysis of the heuristic repair method}\label{sec7}

    The parameter $\Gamma$ represents the number of times the heuristic will be run without improvement before abandoning any further attempts to run the repair heuristic. As mentioned in \cref{sec5}, the parameter $\Gamma$ affects the heuristic repair method's overall performance. This section analyzes the results related to the parameter $\Gamma$.
   
\begin{figure}[htb!]
        \centering
         \includegraphics[scale=0.4]{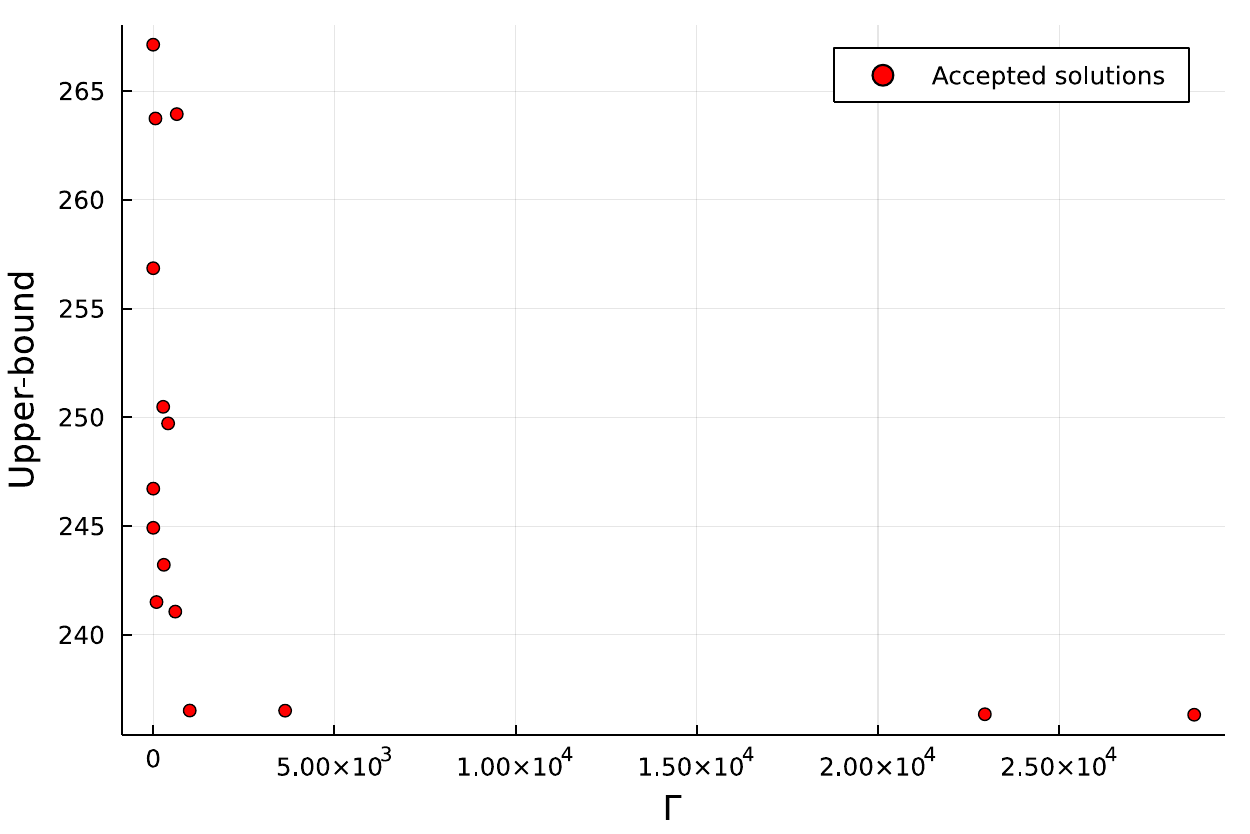}
         \caption{Accepted heuristic repair solutions of instance G:  The x-axis is the number of $\Gamma$ whereas the y-axis is the current upper bound. The red dots represent the heuristic repair solutions. We set a relatively large initial number of $\Gamma = 1 \times 10^6$, and Gurobi accepted 15 heuristic repair solutions. However, the upper bound of the last three solutions decreased slightly with a large $\Gamma$. Further discussion about the selections of $\Gamma$ is worth studying.}
         \label{Fig.7}
    \end{figure}

    \cref{Fig.7} shows the number of accepted heuristic solutions during the branch-and-bound process of Ins.G. It indicates that our heuristic repair method generates 15 solutions that are better than the current solutions found by the solver, which decreases the upper bound quickly. However, \cref{Fig.7} shows the heuristic repair method generates 12 acceptable solutions at $\Gamma < 3000$, while the last three solutions need to set a large $\Gamma$ but return a slightly changed upper bound. Thus, the selection of $\Gamma$ is crucial to the heuristic repair method.

       \cref{parameter test} displays the ``lower bound (LB)'', ``upper bound (UB)'', and ``gap'' between LB and UB for different choices of heuristic repair parameter $\Gamma$ in the same graph instance G. The column ``$\lvert Acc \rvert$" reports the number of accepted heuristic repair solutions.  We use a consistent random number seed = 1.
     
     \begin{table}[ht]
      \centering
        \caption{Ins.G: The column ``LB'' shows the lower bound of instance G, whereas ``UB'' presents the upper bound. The selections of $\Gamma$ between 1000 to 90000. The column ``$\lvert Acc \rvert$" reports the number of accepted heuristic repair solutions, while the ``gap'' displays the deviation between the upper and lower bound.}
        \label{parameter test}
        \begin{tabular}{lllll}
        \hline
        $\Gamma$ & LB       & UB       & gap       & $\lvert Acc \rvert$\\ \hline
        1000       & 199.0254 & 239.5126 & 16.90\% & 11\\
        2000       & 198.7669 & 236.5206 & 15.96\% & 12\\ 
        3000       & 198.8007 & 236.5206 & 15.95\% & 12\\ 
        4000       & 198.7792 & 236.5155 & 15.96\% & 13\\ 
        5000       & 198.7792 & 236.5155 & 15.96\% & 13\\ 
        6000       & 198.6752 & 236.5155 & 16.00\% & 13\\ 
        7000       & 198.7792 & 236.5155 & 15.96\% & 13\\ 
        10000      & 198.6752 & 236.5155 & 16.00\% & 13\\ 
        30000      & 198.4794 & 236.3278 & 16.02\% & 15\\ 
        50000      & 198.4794 & 236.3278 & 16.02\% & 15\\ 
        70000      & 198.4794 & 236.3278 & 16.02\% & 15\\ 
        90000      & 198.4794 & 236.3278 & 16.02\% & 15\\ \hline
        \end{tabular}
    \end{table}

    \cref{parameter test} shows that the best solution is $15.95\%$ at $\Gamma = 3000$ with 12 acceptable heuristic repair solutions. \cref{Fig.7} indicates the first 12 acceptable solutions decrease the upper bound significantly, while only slightly decreasing in the last three solutions. The small number choice of $\Gamma = 1000$ generates 11 acceptable upper bounds, which does not generate a good upper bound for the branch-and-bound tree. Thus, it returns a bad solution of $16.90\%$. Besides, choosing a large number of $\Gamma$ is also expensive in generating unacceptable upper bounds, which are worst than the existing solutions. It is costly to set a large or small number of $\Gamma$, and a selection of 3000 to 4000 seems a reasonable first try for a large-scale problem.
              
%% Conclusion of the paper
\clearpage
\section{Conclusions}\label{sec8}

    In this paper, we proposed a splittable agricultural spraying vehicle routing problem. We are improving the formulation using theoretical insights about the optimal solution structure. This paper provides a SCARP model in two formulations: a basic spray and a large edge demand formulation. The large edge demands model outperforms the basic model in solving problems with high requirements for each edge. Solution methods consist of lazy constraints, symmetries elimination, and heuristic repair. 
    
    Tests on instances generated based on real data show the basic model is limited to small-scale problems, lazy constraints work well on medium-sized problems, and the heuristic repair method outperforms the other methods on large-scale trials. Moreover, the tested results indicated the SCARP model could provide cheaper solutions in some instances compared with the classical CARP. 
    
    The parameter test on $\Gamma$ analyzes parameter selections on the heuristic repair method, which makes the heuristic repair method work efficiently.
    
    Further research can focus on selecting a more precise parameter $\Gamma$ for different problems to improve the efficiency of the heuristic repair method. Another direction is to study the accuracy of the solution to large-scale problems using meta-heuristic methods. The SCARP extension with nursing carts is also worth exploring in the next step.

% \section{Compliance with ethical standards}
%     \begin{itemize}
%         \item Funding: No funding was received for conducting this study.

%         \item Conflict of Interest: Author Qian Wan declares that she has no conflict of interest. Author Rodolfo Garc\'ia-Flores declares that he has no conflict of interest. Author Simon A.Bowly declares that he has no conflict of interest. Author Philip Kilby declares that he has no conflict of interest. Author Andreas T. Ernst declares that he has no conflict of interest.
%         \item Ethical approval: This article does not contain any studies with human participants or animals performed by any of the authors. 
%     \end{itemize}

% Acknowledgements
\section*{Acknowledgments}
I express my profound gratitude to CSIRO Data61 and Monash University for generously awarding me the PhD Scholarship and providing unwavering support throughout my academic journey.

% References
\nocite{*}
% \bibliographystyle{itor}

% Appendix
\vspace*{-10pt}
\appsection{Appendix A}\label{App:A}
    \section{Detailed data tables}\label{secA1}
        \FloatBarrier
        \begin{table}[h]
         \small
         \centering
         \caption{Details of computational results of~\cref{Computational results between Basic model and High models}. In addition to the data mentioned in~\cref{Computational results of different techniques}, this table also provides the number of explored nodes in the column of ``$\lvert ExploredNodes \rvert$'', and the number of simplex iterations in the branch-and-bound process in column ``$\lvert SimplexIter \rvert $''.}
         \label{ Details of computational results between Basic model and High edge-demand model}
           \begin{adjustbox}{max width=\textwidth}
    \begin{tabular}{l l l l}
    \hline
    Ins.   & Reports         & Basic\_model & Large\_model  \\  
    \hline
    \multirow{7}{*}{\begin{tabular}[c]{@{}l@{}}LD\_1\\  sum(D)=260.131 \end{tabular}} & obj & 6110.93 & 6110.93 \\ 
                                 & LB  & 3845.45 & 6052.49\\ 
                                 & r      & 27     & 5        \\ 
                                 & gap    & 37.07\%  & 0.96\%   \\ 
                                 & t      & 7200     & 7006.25   \\  
                                 & $\lvert ExploredNodes \rvert$ & 3432   & 51575  \\ 
                                 & $\lvert SimplexIter \rvert $   & 553273 & 2516396  \\ 
    \hline
    \multirow{7}{*}{\begin{tabular}[c]{@{}l@{}}LD\_2\\  sum(D)=390.195\end{tabular}}
                                 & obj & 13675.18 & 13572.47 \\ 
                                 & LB & 8101.93 & 13503.11\\ 
                                 & r        & 40     & 5        \\  
                                 & gap    & 40.75\%  & 0.51\%   \\ 
                                 & t     & 7200   & 6358.96    \\  
                                 & $\lvert ExploredNodes \rvert$& 1  & 49770  \\ 
                                 & $\lvert SimplexIter \rvert $  & 56616  & 2715323  \\ 
    \hline
    \multirow{7}{*}{\begin{tabular}[c]{@{}l@{}}LD\_3\\  sum(D)=520.26\end{tabular}} 
                                 & obj   & 24058.06  & 23991.43 \\ 
                                 & LB    & 14428.90  & 23832.65 \\ 
                                 & r     &  53  & 4        \\ 
                                 & gap   &40.02 \%  & 0.66\%   \\  
                                 & t     & 7200   & 6014.61 \\ 
                                 & $\lvert ExploredNodes \rvert$ & 79      & 44170   \\ 
                                 & $\lvert SimplexIter \rvert $  & 127326  & 2234217  \\ 
    \hline
    \end{tabular}
    \end{adjustbox}
    \end{table}
\FloatBarrier
\begin{table}
    \centering
    \caption{Details of computational results of~\cref{Computational results of different techniques}. In addition to the data mentioned in~\cref{Computational results of different techniques}, this table also provides the number of explored nodes in the column of ``$\lvert ExploredNodes \rvert$'', and the number of simplex iterations in the branch-and-bound process in column ``$\lvert SimplexIter \rvert $''. }
    \label{Details of computational results of different techniques}
   \begin{adjustbox}{max width=\textwidth}
    \begin{tabular}{l l l l l l l}
    \hline
    Ins. &
        &
        Basic-model &
        Lazy-constraints &
        Symmetries elimination &
        Heuristic repair &
        \multicolumn{1}{c}{H\_no\_sym} \\ \hline
    \multirow{7}{*}{\begin{tabular}[c]{@{}l@{}}Ins.A \\ $\lvert N \rvert$= 16\\ $\lvert E \rvert $=22\\ p=20\\ r=1\\ sum(D)=13.0065\\ no obstacles\end{tabular}} &
        obj &
        32.12 &
        32.12 &
        32.12 &
        32.12 &
        32.12 \\ 
        & LB          & 32.12 & 32.12 & 32.12 & 32.12 & 32.12 \\ 
        & gap             & 0.00\%             & 0.00\%             & 0.00\%             & 0.00\%             & 0.00\%             \\  
        & t(s)               & 0.85               & 2.11               & 1.98        & 5.69               & 4.26               \\ 
        & $\lvert Acc \rvert$             & --                    &  --                  & --                   & 1                  & 1                  \\  
        & $\lvert ExploredNodes \rvert$ & 1                  & 1                  & 1                  & 1                  & 1                  \\  
        & $\lvert SimplexIter \rvert$   & 33                 & 49                 & 37                 & 43                 & 46                 \\ 
    \hline
    \multirow{7}{*}{\begin{tabular}[c]{@{}l@{}}Ins.B \\ $\lvert N \rvert$ =20\\ $\lvert E \rvert$=28\\ p=20\\ r=1\\ sum(D)=17.139\\ no obstacles\end{tabular}} &
        obj &
        42.64 &
        42.64 &
        42.64 &
        42.64 &
        42.64 \\ 
        & LB          & 42.64 & 42.64 & 42.64 & 42.64 & 42.64 \\ 
        & gap             & 0.00\%             & 0.00\%             & 0.00\%             & 0.00\%             & 0.00\%             \\  
        & t(s)               & 28.01              & 2.07               & 2.04               & 5.97               & 4.03               \\  
        & $\lvert Acc \rvert$             & --                   & --                   & --                   & 1                  & 1                  \\  
        &  $\lvert ExploredNodes \rvert$  & 1                  & 1                  & 7                  & 5                  & 1                  \\  
        & $\lvert SimplexIter \rvert$   & 71                 & 68                 & 79                 & 124                & 82                 \\ 
    \hline
    \multirow{7}{*}{\begin{tabular}[c]{@{}l@{}}Ins.C\\ $\lvert N \rvert$ =30\\ $\lvert E \rvert$=43\\ p=20\\ r=2\\ sum(D)=28.362\\ no obstacles\end{tabular}} &
        obj &
        fail &
        76.15 &
        76.15 &
        76.15 &
        76.15 \\  
        & LB          & --                  & 76.15 & 76.15 & 76.15 & 76.15 \\ 
        & gap             & --                  & 0.00\%             & 0.00\%             & 0.00\%             & 0.00\%             \\ 
        & t(s)               & --                  & 49.99       & 62.07              & 79.04              & 91.90       \\ 
        & $\lvert Acc \rvert$             & --                  & --                   & --                   & 3.33        & 3                  \\  
        & $\lvert ExploredNodes \rvert$ & --                  & 317766             & 375578             & 420542             & 535452             \\ 
        & $\lvert SimplexIter \rvert$   & --                  & 3977670            & 4971315            & 5590590            & 6569088            \\ 
    \hline
    \multirow{7}{*}{\begin{tabular}[c]{@{}l@{}}Ins.D\\ $\lvert N \rvert$ =42\\ $\lvert E \rvert$=61\\ p=20\\ r=2\\ sum(D)=32.5815\\ with a obstacle\end{tabular}} &
        obj &
        fail &
        84.30 &
        84.30 &
        84.30 &
        84.30 \\ 
        & LB          & --                  & 84.30 & 84.30 & 84.30 & 84.30 \\ 
        & gap             & --                  & 0.00\%             & 0.00\%             & 0.00\%             & 0.00\%             \\ 
        & t(s)              & --                  & 5.74        & 5.62               & 10.93              & 9.23               \\  
        & $\lvert Acc \rvert$             & --                  & --                  & --                   & 6                  & 7                  \\  
        & $\lvert ExploredNodes \rvert$ & --                  & 5761               & 4424               & 7742               & 7656               \\ 
        & $\lvert SimplexIter \rvert$   & --                  & 143823             & 109771             & 168798             & 178461             \\ 
    \hline
    \multirow{7}{*}{\begin{tabular}[c]{@{}l@{}}Ins.E\\ $\lvert N \rvert$ =56\\ $\lvert E \rvert$=82\\ p=20\\ r=3\\ sum(D)=46.59\\ with a obstacle\end{tabular}} &
        obj &
        fail &
        127.65 &
        127.65 &
        127.65 &
        127.65 \\ 
        & LB          & --                  & 124.23 & 123.03 & 124.52 & 124.14 \\  
        & gap             & --                  & 2.68\%             & 3.62\%             & 2.45\%             & 2.75\%             \\  
        & t(s)               & --                  & 7200               & 7200               & 7200               & 7200               \\  
        & $\lvert Acc \rvert$             & --                  &--                    & --                   & 8                  & 6                  \\  
        & $\lvert ExploredNodes \rvert$ & --                  & 4509341            & 4298279            & 5212149            & 5145643            \\ 
        & $\lvert SimplexIter \rvert$   & --                  & 371701414          & 349635575          & 447436137          & 426000722          \\ 
    \hline
    \multirow{7}{*}{\begin{tabular}[c]{@{}l@{}}Ins.F\\ $\lvert N \rvert$ =78\\ $\lvert E \rvert$=115\\ p=20\\ r=3\\ sum(D)=52.983\\ with two obstacles\end{tabular}} &
        obj &
        fail &
        150.74 &
        151.98 &
        152.19 &
        149.79 \\ 
        & LB          & --                  & 140.31 & 140.67 & 140.57 & 140.69 \\ 
        & gap             & --                  & 6.91\%             & 7.42\%             & 7.64\%             & 6.06\%             \\ 
        & t(s)               & --                  & 7200               & 7200               & 7200               & 7200               \\ 
        & $\lvert Acc \rvert$             & --                  &   --                 &   --                 & 8        & 8                  \\  
        & $\lvert ExploredNodes \rvert$ & --                  & 322947             & 425121             & 444244             & 400689             \\  
        & $\lvert SimplexIter \rvert$   & --                  & 47437951           & 63547502           & 61216769           & 63985627           \\ 
    \hline
    \multirow{7}{*}{\begin{tabular}[c]{@{}l@{}}Ins.G\\ $\lvert N \rvert$ =98\\ $\lvert E \rvert$=145\\ p=20\\ r=4\\ sum(D)=74.7765\\ with two obstacles\end{tabular}} &
        obj &
        fail &
        fail &
        fail &
        235.77 &
        251.58 \\ 
        & LB          & --                  & --                  & --                  & 198.60 & 199.21 \\  
        & gap             & --                  & --                 & --                  & 15.75\%            & 20.72\%            \\  
        & t(s)               & --                 & --                  & --                  & 7200               & 7200               \\ 
        & $\lvert Acc \rvert$             & --                 & --                  & --                  & 13                 & 6                  \\ 
        & $\lvert ExploredNodes \rvert$ & --                  & --                  & --                  & 177727             & 175286             \\ 
        & $\lvert SimplexIter \rvert$   & --                  & --                  & --                  & 20497276           & 17737926           \\ 
    \hline
    \end{tabular}
    \end{adjustbox}
    \end{table}  

\end{document}